\documentclass[a4paper,12pt,reqno]{amsart}
\usepackage{graphics}
\newtheorem{lem}{Lemma}
\newtheorem{theo}{Theorem}
\newtheorem{rem}{Remark}
\newtheorem{cor}{Corollary}
\newtheorem{defi}{Definition}
\numberwithin{equation}{section}

\begin{document}

\title[Even circuits]{Even circuits of prescribed clockwise parity}

\author[Ilse Fischer and C.H.C. Little]{\box\Adr}

\newbox\Adr
\setbox\Adr\vbox{
\centerline{ Ilse Fischer}
\centerline{Universit\"at Klagenfurt, A-9020 Klagenfurt, Austria}
\centerline{{\tt Ilse.Fischer@uni-klu.ac.at}}
\vspace{0.2cm}
\centerline{and}
\vspace{0.2cm}
\centerline{ C.H.C. Little}
\centerline{Massey University, Palmerston North, New Zealand}
\centerline{{\tt C.Little@massey.ac.nz}}
}

\maketitle

\begin{abstract}
We show that a graph has an orientation under which every circuit
of even length is clockwise odd if and only if the graph contains
no subgraph which is, after the contraction of at most one
circuit of odd length, an even subdivision of $K_{2,3}$. In fact we
give a more general
characterisation of graphs that have an orientation under
which every even circuit has a prescribed clockwise parity.
This problem was
motivated by the study of {\it Pfaffian} graphs, which are the
graphs that have an orientation under which every alternating
circuit is clockwise odd. Their significance is that they are
precisely the graphs to which Kasteleyn's powerful method \cite{kasteleyn}
for enumerating perfect matchings may be applied.
\end{abstract}

\section{Introduction}

Consider the three (even) circuits in $K_{2,3}$. Is it possible to
find an orientation under which all these circuits are clockwise
odd, if the clockwise parity of a circuit of even length is defined as the
parity of the number of edges that are directed in agreement with
a specified sense? However $K_{2,3}$ is oriented
one observes that the total number of clockwise even circuits is odd
and therefore it is not possible to find such an orientation. In
this paper we present a characterisation, in terms of forbidden
subgraphs,
of the graphs that have an orientation under which every even circuit is
clockwise odd. It will turn out that the
non-existence of such an orientation can in a sense always be put down to an
even subdivision of $K_{2,3}$. (See Corollary~\ref{even}.)

\medskip

We were motivated to study this problem by our work on a
characterisation of Pfaffian graphs. A {\it Pfaffian
orientation} of a graph is an orientation under which every
alternating circuit is clockwise odd, an alternating circuit being
a circuit which is the symmetric difference of two
perfect matchings. A {\it Pfaffian} graph is a graph that admits a
Pfaffian orientation. In \cite{kasteleyn} Kasteleyn introduced a
remarkable method for enumerating perfect matchings in Pfaffian
graphs, reducing the enumeration to the evaluation of the
determinant of the skew adjacency matrix of the Pfaffian directed graph.
He has shown that all planar graphs are Pfaffian. However
a general characterisation of Pfaffian graphs is still not known. For
research in this direction see
\cite{Li1,Li2,LiReFi,FiLi,RoSeTh}.

\medskip

Our characterisation of the graphs that admit an orientation under
which every even circuit is clockwise odd will be an easy
consequence of our main theorem, which gives a more general 
characterisation of the
graphs that
have an orientation under which every even circuit has a
prescribed (not necessarily odd) clockwise parity.
Before we are able to state our main theorem
we need some definitions.

\begin{defi}
Let $G$ be a graph and $J$ an assignment of clockwise
parities to the even circuits of $G$. An orientation of $G$ is said to be
$J$-compatible if every even circuit of $G$ has the clockwise parity
prescribed by $J$. Otherwise the orientation is $J$-incompatible. The 
graph $G$ is
said to be $J$-compatible if $G$ admits a $J$-compatible orientation,
and $J$-incompatible otherwise.
\end{defi}

Our main theorem (Theorem~\ref{main}) characterises $J$-compatible 
graphs in terms of
forbidden subgraphs. Before we are able to formulate it, we have
to introduce two relevant graph operations. To this end we need
the following fundamental definition and lemma.

\begin{defi}
\label{intractabledefinition}
Let $G$ be a graph and $J$ an assignment of clockwise parities to
the even circuits of $G$. A set $\mathcal S$ of even circuits in $G$ is said to
be {\it $J$-intractable} if the symmetric difference of the
circuits in $\mathcal S$ is empty and the parity of the number of
clockwise even circuits in $\mathcal S$ with respect to an 
orientation is unequal
to the parity of the number of clockwise even circuits
in $\mathcal S$ with respect to the assignment $J$.
\end{defi}

Observe that the parity of the number of clockwise even circuits with respect
to an orientation in a
$J$-intractable set ${\mathcal S}$ does not depend on the
orientation since the reorientation of a single edge changes the
clockwise parity of an even number of circuits in ${\mathcal S}$.

\begin{lem}
\label{intractable}
  Let $G$ be a graph and $J$ an assignment of clockwise
parities to the even circuits of $G$. Then $G$ is $J$-incompatible if
and only if $G$ contains a $J$-intractable set of even circuits.
\end{lem}

{\it Proof.} The fact that the existence of a $J$-intractable set
implies that $G$ is $J$-incompatible follows from the remark above 
the formulation
of the lemma.

Suppose that $G$ is $J$-incompatible and orient $G$ arbitrarily. The
existence of a $J$-compatible orientation of $G$ is equivalent to the
solvability of a certain system of linear equations over the field 
$\mathbb{F}_2$.
In these equations the variables
correspond to the edges of the even circuits of $G$. For every even
circuit $C$
there is a corresponding equation in which the sum of the variables
corresponding to the edges of $C$ is $1$ if
and only if the clockwise parity of $C$ is not that
prescribed by $J$. A solution of this system is an assignment of
zeros and ones to the edges of the even circuits of $G$. A
$J$-compatible orientation of
$G$ can be obtained from the fixed orientation by reorienting
precisely those edges to which the solution assigns a $1$. The
lemma now follows from the solvability criteria for
systems of linear equations. \qed

\medskip

\begin{figure}
\begin{center}
\setlength{\unitlength}{1cm}
\hspace{0.5cm}\mbox{\scalebox{0.50}{%
  \includegraphics{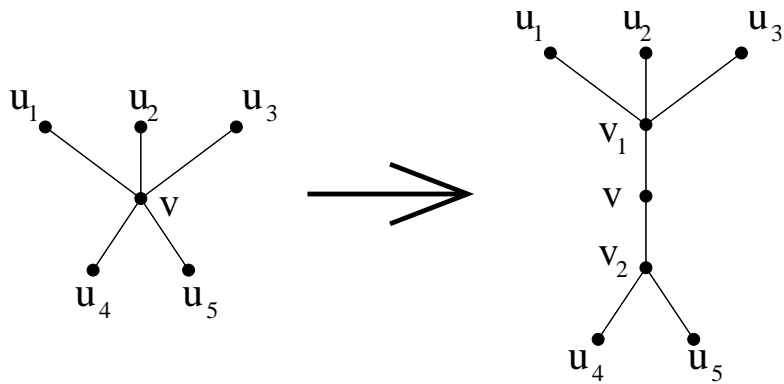}}}%
\end{center}
\caption{Split vertices to obtain an even vertex splitting.}
\label{split}
\end{figure}

Let $G$ be a graph and let $H$ be a graph obtained from $G$ by the
contraction of the two edges $e$ and $f$ incident on some vertex $v$
in $G$ of degree $2$. Thus $EH = EG - \{e,f\}$. We may describe $G$ as
an {\it even vertex splitting} of $H$. (See Figure~\ref{split}.) Any 
even circuit $C_H$ in
$H$ is the intersection with $EH$ of a unique even circuit $C$ in $G$. To
any assignment $J$ of clockwise parities to the even circuits in $G$
there corresponds an assignment $J_H$ of clockwise parities to the
even circuits in $H$ so that any even circuit $C_H$ in $H$ is assigned the
same clockwise parity as $C$ in $G$. We say that $J_H$ is {\it induced}
by $J$. If either $e$ or $f$ is incident on a vertex of degree $2$
other than $v$, then it is also true that the intersection with $EH$
of any even circuit $C$ in $G$ yields an even circuit $C_H$ in $H$. In
this case any assignment $J_H$ of clockwise parities to the even
circuits in $H$ corresponds to a unique assignment $J$ of clockwise parities
to the even circuits in $G$ so that $J_H$ is the assignment induced
by $J$. We then say that $J$ is also {\it induced} by $J_H$.

Similarly let $H$ be obtained from $G$ by contracting a circuit $A$ of odd
length. Thus $EH = EG - A$. Any even circuit $C_H$ in $H$ is the intersection
with $EH$ of a unique even circuit $C$ in $G$: we have
$C \cap EH = C_H$ and if $C \ne C_H$ then $C \cap A$ is the path of 
even length in $A$
joining the ends of the path $C_H$ in $G$. To
any assignment $J$ of clockwise parities to the even circuits in $G$
there corresponds an assignment $J_H$ of clockwise parities to the
even circuits in $H$ so that any even circuit $C_H$ in $H$ is assigned the
same clockwise parity as $C$ in $G$. We say that $J_H$ is {\it induced}
by $J$.

In the following lemma we summarise some basic facts:

\begin{lem}
\label{basic}
  Let $G$ be a graph and $J$ an assignment of clockwise parities to
the even circuits of $G$.

(1) Let $H$ be a subgraph of $G$ and $J_H$ the restriction of $J$
to the even circuits of $H$. If $G$ is $J$-compatible then $H$ is
$J_H$-compatible.

(2) Let $H$ be obtained from $G$ by contracting the two edges
incident on a vertex of degree $2$.
The assignment $J$ induces an assignment $J_H$ of
clockwise parities to the even circuits in $H$. If $G$ is
$J$-compatible
then $H$ is $J_{H}$-compatible. If either of the two contracted edges is
incident on another vertex of degree $2$ then $G$ is $J$-compatible if and
only if $H$ is $J_H$-compatible.

(3) Let $H$ be obtained from $G$ by contracting a circuit of odd length. The
assignment $J$ induces an assignment $J_H$ of
clockwise parities to the even circuits in $H$. If $G$ is
$J$-compatible then $H$ is
$J_H$-compatible.
\end{lem}

{\it Proof.} (2) Every $J_H$-intractable set of even circuits in
$H$ induces a $J$-intractable set of even circuits in $G$. If either
of the two contracted edges is incident on another vertex of
degree two then every $J$-intractable set of even circuits in $G$ also
induces a $J_H$-intractable set of even circuits in $H$.

\smallskip

(3) Every $J_H$-intractable set of even circuits in $H$ induces
a $J$-intractable set of even circuits in $G$. (Note that if the
symmetric difference of some even circuits in $H$ is empty, then
the symmetric difference of the corresponding even circuits in $G$
is empty as well, for it is obvious that this symmetric difference
is both an even cycle and a subset of the odd circuit $EG - EH$.) \qed

\medskip

In the following three paragraphs we introduce the minimal
$J$-incompatible
graphs which we need in the formulation of our main theorem. We say
that an assignment $J$ is {\it odd} or {\it even} if it assigns,
respectively, an odd or an even clockwise parity to every even
circuit.

\smallskip

Let $O_1=K_{2,3}$ and let $O_2$ be the graph we obtain from $K_4$ by
subdividing once all edges incident on one fixed vertex. (See
Figure~\ref{g1}.) Observe that
$O_1$ and $O_2$ are $J$-incompatible with respect to the odd assignment $J$.
In fact $O_{1}$ and $O_{2}$ are $J$-incompatible precisely
for those assignments $J$ that prescribe an even number of even
circuits of these graphs to be of even clockwise parity.
For these assignments Lemma~\ref{basic}(3) shows
that the $J$-incompatibility of $O_2$
can be attributed to the fact that $O_1$ is
$J$-incompatible, since the contraction of the triangle in $O_2$ gives
$O_1$.

\smallskip

Let $E_1$ be the graph consisting of two vertices and three edges
joining them, let $E_2=K_4$ and let $E_3$ be the graph we obtain from
$K_4$ by subdividing each edge in a fixed even circuit once. (See
Figure~\ref{g1}.) Then
$E_1$, $E_2$ and $E_3$ are $J$-incompatible with respect to the even
assignment $J$. More generally $E_{1}$, $E_{2}$ and $E_{3}$ are
$J$-incompatible precisely for those assignments $J$ that prescribe an
odd
number of even circuits to be of even clockwise parity.
Again by Lemma~\ref{basic}(3) the fact that $E_2$
is $J$-incompatible can be put down to the fact that $E_1$ is
$J$-incompatible, since the contraction of a triangle in $E_2$ gives
$E_1$.

\smallskip

\begin{figure}
\begin{center}
\setlength{\unitlength}{1cm}
\hspace{0.5cm}\mbox{\scalebox{0.30}{
  \includegraphics{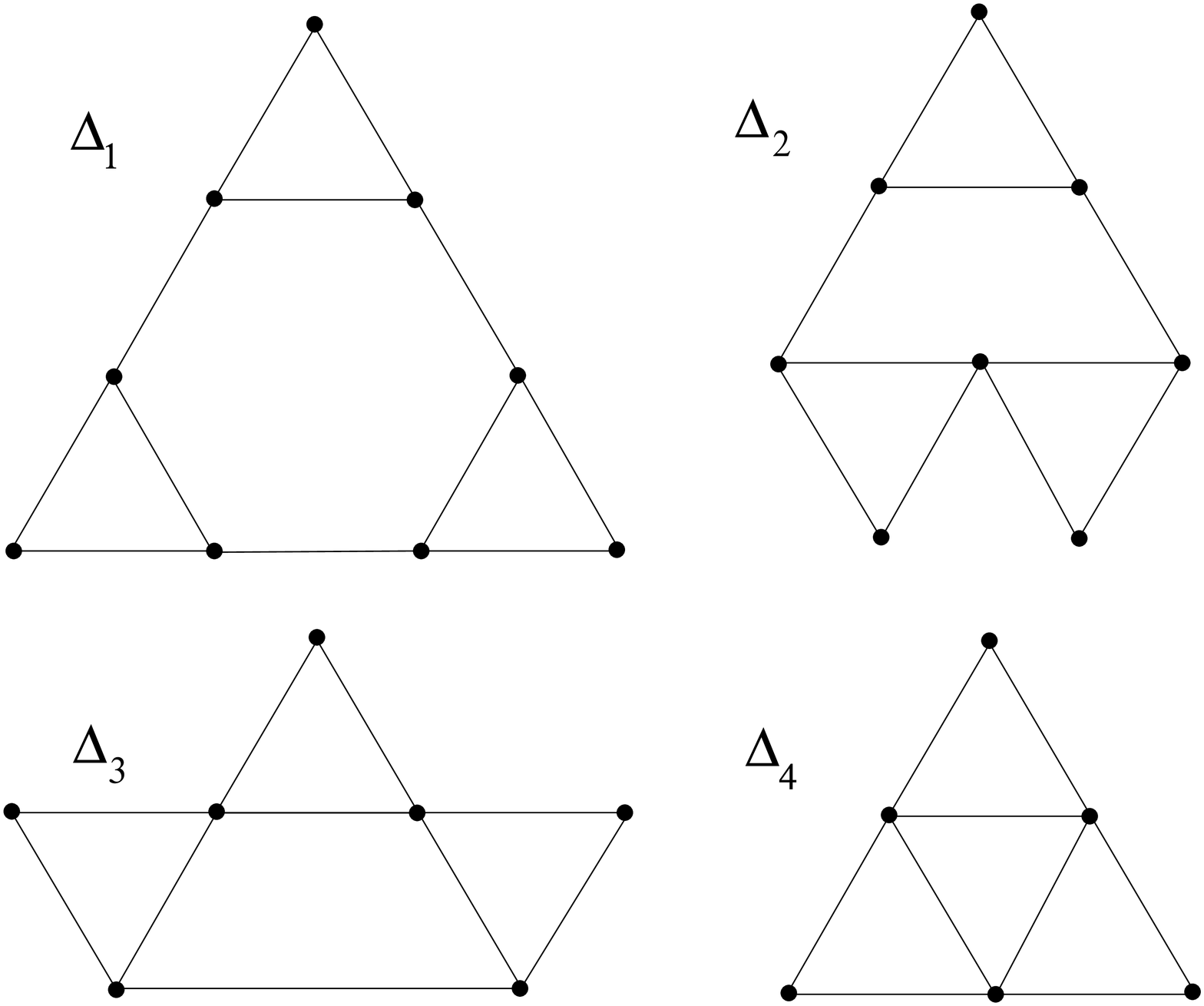}}}%
\end{center}
\caption{The graphs $\Delta_1, \Delta_2, \Delta_3, \Delta_4$.}
\label{3t}
\end{figure}

Consider the four graphs $\Delta_1, \Delta_2, \Delta_3, \Delta_4$ in
Figure~\ref{3t}. Note that the last three are obtained from
$\Delta_1$ by contracting edges. Each of these four graphs
has exactly four even circuits and is $J$-incompatible if and
only if $J$ prescribes an odd number of them to be clockwise even.
This observation follows from Lemma~\ref{intractable} because in each of
these graphs
the set of all even circuits is the only non-empty dependent set of 
even circuits
with respect to symmetric difference.

\medskip

Let $G$ be a graph and let $G_0, G_1, \ldots, G_k$ be graphs such
that $G_0 = G$ and, for each $i > 0$, the graph $G_i$ is an even
vertex splitting of $G_{i-1}$. Then $G_k$ is said to be an {\it even
splitting} of $G$. There is a special case in which, for each
$i > 0$, $G_i$ can be obtained from $G_{i-1}$ by subdividing an edge
twice. In this case we describe $G_k$ as an {\it even subdivision}
of $G$. If no vertex of $G$ is of degree greater than $3$, then each
even splitting of $G$ is also an even subdivision of $G$. If $H$ is 
an even splitting of
$G$ and $J_H$ is an assignment
of clockwise parities to the even circuits in $H$, then we may apply
the definition of an induced assignment inductively to obtain an
assignment $J$ of clockwise parities to the even circuits in $G$. This
assignment is also said to be {\it induced} by $J_H$. By applying 
Lemma~\ref{basic}(2)
inductively we find that if $H$ is $J_H$-compatible then $G$ is
$J$-compatible. Thus if $G$ is $J$-incompatible then $H$ is
$J_H$-incompatible. The converse also holds if $H$ is an even
subdivision of $G$.

\medskip

Now we formulate our main theorem.

\begin{theo}
\label{main} Let $G$ be a graph and $J$ an assignment of clockwise
parities to the even circuits of $G$. Then $G$ is $J$-incompatible if
and only if $G$ contains an even splitting $H$ of one of $O_1$,
$O_2$, $E_1$, $E_2$, $E_3$, $\Delta_1$, $\Delta_2$, $\Delta_3$, 
$\Delta_4$ with the following
property: if $H$ is an even splitting of $O_i$ for some $i$ then $J$ 
prescribes an even
number of clockwise even circuits to the three even circuits of $H$, if
$H$ is an even splitting of $E_i$ for some $i$ then $J$ prescribes an
odd number of clockwise
even circuits to the three even circuits of $H$ and if $H$ is an even splitting
of $\Delta_i$ for some $i$ then $J$ prescribes an odd number of 
clockwise even circuits
to the four even circuits of $H$.
\end{theo}

The ``if'' direction in the theorem is obvious by Lemma~\ref{intractable} and
Lemma~\ref{basic}. We obtain the following immediate corollaries.

\begin{cor}
\label{even}
A graph $G$ does not admit an orientation under which every even
circuit is clockwise odd if and only if it contains a subgraph which
is, after the contraction of at most one odd circuit, an even
subdivision of $K_{2,3}$.
\end{cor}

{\it Proof.} For each even splitting of $E_1$, $E_2$,
$E_3$, $\Delta_1$, $\Delta_2$, $\Delta_3$ or $\Delta_4$ contained
in $G$ the odd
assignment prescribes an even number of clockwise even circuits to
its set of even circuits. \qed

\begin{cor}
A graph $G$ does not admit an orientation under which every even
circuit is clockwise even if and only if it contains a subgraph
which is, after the contraction of at most one odd circuit,
an even subdivision of $E_1$ or $E_3$.
\end{cor}

{\it Proof.} For each even splitting of $O_1$ or $O_2$
contained in $G$ the even assignment prescribes an odd number of
clockwise even circuits to its set of even circuits, for both graphs
have three even circuits. Moreover for each even
splitting of $\Delta_1$, $\Delta_2$, $\Delta_3$ or $\Delta_4$
contained in $G$ the even assignment
prescribes an even number of clockwise even circuits to its set
of even circuits, for these graphs each have four even circuits. \qed

\section{An arc decomposition theorem}

Circuits, non-empty paths and, more generally, subgraphs without
isolated vertices are determined by their edge sets and are
therefore identified with them in this paper. If $u$ and $v$ are
vertices of a path $P$, then $P[u,v]$ denotes the subpath of $P$ that
joins $u$ and $v$. If $G$ is a graph
and $V'$ is a subset of the vertex set $VG$ of $G$ then $G[V']$
denotes the subgraph of $G$ spanned by $V'$.
Similarly  if $E'$ is a subset of the edge set $EG$ of $G$ then
$G[E']$ denotes the subgraph of $G$ spanned by $E'$.

Let $H_1$ and $H_2$ be two sets of edges in $G$. An $H_1
\overline{H_2}$-arc is a path in $H_1$ which joins two distinct
vertices in $VG[H_2]$ but does not have an inner vertex in $VG[H_2]$. A
$G \overline{H_2}$-arc is also called an $\overline{H_2}$-arc.

\begin{defi}
A graph $G$ without isolated vertices is said to be {\it even-circuit-connected}
if for every bipartition $\{H_{1},H_{2}\}$ of $EG$ there
exists an even circuit $C$ which meets $H_{1}$ and $H_{2}$.
\end{defi}

Every even-circuit-connected graph is $2$-connected. Indeed, suppose 
there exists a
vertex $v$ such that $G-\{v\}$ is disconnected. Let
$H_{1},H_{2},\dots,H_{k}$ be the components of $G-\{v\}$ and let
$H'_{i} = G[VH_{i} \cup \{v\}]$ for each $i$. Then for every circuit $C$
of $G$ there exists an $i$, $1 \le i \le k$, with $C \subseteq
EH_{i}'$, a contradiction.

\medskip

First we prove a decomposition theorem on even-circuit-connected 
graphs. Note that in
our characterisation of $J$-incompatible graphs in terms of forbidden
subgraphs, even-circuit-connected graphs are the only graphs of interest
since every $J$-incompatible graph that is minimal with respect to 
edge deletion
is clearly even-circuit-connected.

\smallskip

Let $H$ be a subgraph of $G$ and $C$ an even circuit in $G$ which includes
$EG - EH$ and meets $EH$. If there are $n$ $C\overline{H}$-arcs, then $G$ is
said to be obtained from $H$ by an {\it $n$-arc adjunction}. An
{\it arc decomposition} of an even-circuit-connected
graph $G$ is a sequence $G_{0}, G_{1}, \dots, G_{k}$ of 
even-circuit-connected subgraphs of
$G$ such that $EG_{0}$ is an even circuit, $G_{k} = G$ and, for every $i > 0$,
$G_{i}$ is obtained from $G_{i-1}$ by an $n$-arc adjunction with
$n=1$ or $n=2$. Moreover we assume that, for each $i$, every even circuit in
$G_{i}$ which meets $EG_{i} - EG_{i-1}$
contains $EG_{i} - EG_{i-1}$. We shall show that every 
even-circuit-connected graph has an
arc decomposition. For this purpose we need the following version of 
Menger's theorem.

\begin{theo}
\label{menger}
\cite{Ha} Let $S$ and $T$ be sets of at least $n$ vertices in an $n$-connected
graph $G$. Then there are $n$ vertex disjoint paths joining vertices in $S$
to vertices in $T$ such that no inner vertex of these paths is in
$VS \cup VT$.
\end{theo}

\begin{lem}
\label{decomp} Let $H$ be a non-empty proper even-circuit-connected 
subgraph of an
even-circuit-connected graph $G$. Then $G$ has an even circuit $C$ 
that meets $EH$, admits
just one or two $C \overline{H}$-arcs and has the property that $G[H 
\cup C]$ is
even-circuit-connected.
Moreover, if $G$ is bipartite or $H$ is not, then $C$
may be chosen to admit just  one  $C \overline{H}$-arc.
\end{lem}

{\it Proof.} Suppose first that $G$ is bipartite. By hypothesis
there is an edge $e$ in $EG - EH$.  By the $2$-connectedness of $G$
and Theorem~\ref{menger} there are vertex disjoint paths $P$ and $Q$
in $EG - EH$ joining the ends of $e$ to two distinct vertices $u$ and
$v$, respectively, in $VH$ such that neither $P$ nor $Q$ has an inner vertex in
$VH$. Since $H$ is even-circuit-connected and therefore connected, a 
path $R$ in $H$ joins
$u$ to $v$. Thus $P \cup \{e\} \cup Q \cup R$ is a circuit $C$ in $G$
meeting $EH$ (since $u \not= v$) and having $P \cup \{e\} \cup Q$ as its unique
$C \overline{H}$-arc. Moreover $C$ is even since $G$ is bipartite.

Suppose therefore that $G$ is not bipartite. Again we may
construct the circuit $C$ as in the previous case, and the proof
is complete if $C$ is even. Suppose therefore that $C$ cannot be
chosen to be even. Since $G$ is even-circuit-connected there exists 
an even circuit
$D$ which meets $EH$ and $EG - EH$. Let $S$ and $T$ be two
distinct $D \overline{H}$-arcs, joining $w$ to $x$ and $y$ to $z$,
respectively. The fact that $D$ meets $EH$ implies $w \not=
x$, $y \not=z$ and $\{w, x\} \not= \{y,z\}$. Let $U$ be a path in
$H$ joining $w$ to $x$. Since $H$ is 2-connected there exist two
vertex disjoint paths $V$ and $W$ in $H$ joining $y$ and $z$, respectively, 
to distinct
vertices of $U$ and such that neither has an inner vertex
in $VU$. Let $s$ and $t$ be the ends of $V$ and $W$, respectively, in $VU$. By
assumption $S \cup U$ is an odd circuit  and therefore it
includes a path $X$, joining $s$ to $t$, such that
$$
|X| \equiv |T| + |V| + |W| \qquad (\text{mod} \hspace{2mm} 2).
$$
Then $T \cup V \cup W \cup X$ is a circuit $C$ of even length, and the only
$C\overline{H}$-arcs are $T$ and possibly $S$.

Finally suppose that $H$ is not bipartite. Therefore $H$ has an odd
circuit $O$. Since $H$ is even-circuit-connected and therefore $2$-connected,
there are vertex disjoint paths
$M$ and $N$ in $H$ joining $w$ and $x$, respectively, to distinct
vertices $p$ and $q$ in $VO$ but having no inner
vertex in $VO$. Since $O$ is odd, it includes a path $Y$ joining $p$
and $q$ such that
$$
|Y| \equiv |M| + |N| + |S| \qquad (\text{mod} \hspace{2mm} 2).
$$
Then $M \cup N \cup S \cup Y$ is an even circuit $C$ in $G$, and
$S$ is the unique $C \overline{H}$-arc.

It remains to show that $G[H \cup C]=:H'$ is even-circuit-connected. 
Suppose that
$\{K_{1},K_{2}\}$
is a bipartition of $EH'$. If $K_{1} \supseteq EH$ and $K_{2} 
\subseteq EH' - EH$
then $C$ is an even circuit which meets $K_{1}$ and $K_{2}$. Thus we 
may assume that
$K_{l} \cap EH =:K'_{l} \not= \emptyset$ for $l=1,2$. By the assumption that
$H$ is even-circuit-connected there exists an even circuit in $H$ 
which meets $K'_{1}$ and $K'_{2}$, and
therefore $K_{1}$ and $K_{2}$. \qed

\medskip

Lemma~\ref{decomp} shows that every even-circuit-connected graph $G$ 
has an arc decomposition
$G_{0},G_{1},\dots,G_{n}$ with at most one $2$-arc adjunction. The single
$2$-arc adjunction is necessary if and only if $G$ is not bipartite. In
this case the arc decomposition can be chosen so that $G_{1}$ is obtained
from $G_{0}$ by a $2$-arc adjunction as we see in the following theorem.

\begin{theo}
\label{optimal}
An even-circuit-connected graph $G$ has an arc decomposition $G_{0}, 
G_{1}, \dots, G_{k}$ such
that $G_{i}$ is obtained from $G_{i-1}$ by a single arc adjunction for all
$i>1$.
\end{theo}

{\it Proof.} By Lemma~\ref{decomp} let $H_{0}, H_{1}, \dots, H_{n}$ be an arc
decomposition of $G$
such that $H_{i}$ is obtained from $H_{i-1}$ by a 2-arc adjunction for some
$i>1$ and $H_{j}$ is obtained from $H_{j-1}$ by a single arc adjunction for
all $j \not= i$. Let $EH_{i} - EH_{i-1} = P \cup Q$, where $P$ and $Q$ are the
two $H_{i} \overline{H_{{i-1}}}$-arcs. Let $P$ join $w$ to $x$ and 
$Q$ join $y$ to $z$. We
distinguish cases according to whether $w,x,y,z$ are distinct.

{\it Case~1.} Suppose that $w,x,y,z$ are distinct. Since $H_{i-1}$ is 
2-connected,
we may assume by
Theorem~\ref{menger}, the symmetry of $w$ and $x$ and the symmetry of $y$ and
$z$ that $H_{i-1}$
has vertex disjoint paths $R$ joining $w$ to $y$ and $S$ joining $x$ 
to $z$. Similarly there
are two vertex disjoint paths $T$ and $U$ in $H_{i-1}$ joining vertices in
$VR$ to vertices in $VS$ but having no inner vertex in
$VR \cup VS$. Let $T$ join vertex $a$ in $VR$ to vertex $b$ in $VS$ and
let $U$ join vertex $c$ in $VR$ to vertex $d$ in $VS$. With no less generality
we may assume that $c \in VR[a,y]$. Then $R[c,a] \cup T \cup S[b,d] \cup U$ is
an even circuit $C$, since $H_{i-1}$ is bipartite. Note also
that $P \cup R[w,a] \cup T \cup S[b,x]$ and $Q \cup R[y,c] \cup U \cup S[d,z]$
are odd circuits $A$ and $B$, respectively, for neither $G[H_{i-1} \cup P]$
nor $G[H_{i-1} \cup Q]$ is an even-circuit-connected graph.

Let $G_{0}: = G[C]$  and $D = P \cup R \cup Q \cup
S$. Then $D$ is an even circuit since $D = A + B + C$.
Furthermore observe that $G_{1}: = G[C \cup D]$ is a non-bipartite 
even-circuit-connected
graph obtained from $G_{0}$ by a $2$-arc adjunction.
Thus the assertion follows from Lemma~\ref{decomp}.

{\it Case~2.} In the remaining case observe that $w \not= x$, $y \not= z$ and
$\{w,x\} \not= \{y,z\}$ for there exists an even circuit which includes $P
\cup Q$ and meets $EH_{i-1}$. Thus we may assume that $x=y$
and $|\{w,x,z\}| = 3$ without loss  of generality.
If there are edges $e$ and $f$ in $EH_{i-1}$ joining $x$ to $w$ and
$z$ respectively, then
set $R=\{e\}$ and $S=\{f\}$. Otherwise,
since $H_{i-1}$ is 2-connected, $x$ is joined in $H_{i-1}$ by an
edge $g$ to a vertex $v$
not in $\{w,x,z\}$. Without loss of generality we assume that there are
vertex disjoint paths $R'$ and $S$ joining $v$ to $w$ and $x$ to $z$, 
respectively.
Set $R= R' \cup \{g\}$. By the $2$-connectedness of $H_{i-1}$
there exists a path $T$ in $H_{i-1} - \{x\}$ joining a vertex $a$ in $VR$
to a vertex $b$ in $VS$ but having no
inner vertex in $VR \cup VS$. Then $C:=R[a,x] \cup S[x,b] \cup T$
is an even circuit for $H_{i-1}$ is bipartite. Set
$$
D=P \cup R[w,a] \cup T \cup S[b,z] \cup Q.
$$
Observe that $D$ is an even circuit since $D= C + P + R  + Q + S$
and $P+R$ and $Q + S$ are odd circuits. Finally set $G_{0} = G [C]$ and
$G_{1} = G[C \cup D]$. Then $G_{1}$ is a non-bipartite 
even-circuit-connected graph which can be
obtained from $G_{0}$ by a 2-arc adjunction. Again the assertion follows
from Lemma~\ref{decomp}. \qed

\medskip

\begin{figure}
\begin{center}
\setlength{\unitlength}{1cm}
\hspace{0.5cm}\mbox{\scalebox{0.50}{%
  \includegraphics{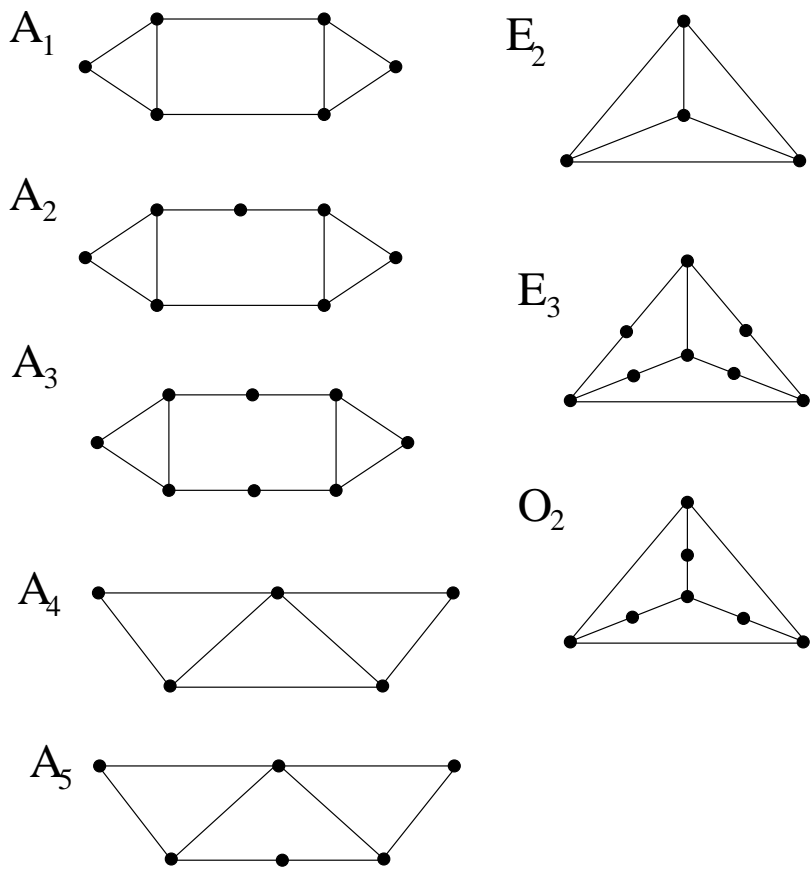}}}%
\end{center}
\caption{$G_{1}$ in Theorem~\ref{optimal} is an even subdivision of 
one of these
   graphs.}
\label{g1}
\end{figure}

\begin{rem}
In the previous proof $G_{1}$ is an even subdivision of one of the
graphs in Figure~\ref{g1}.
\end{rem}

\section{Proof of Theorem~\ref{main}}

For the rest of the paper let $G$ be a graph and $J$ an assignment of 
clockwise parities to
the even circuits of $G$. Assume that $G$ is minimally 
$J$-incompatible with respect
to the deletion of an edge. Let $G_{0}, G_{1}, \dots, G_{k}$ be an arc
decomposition of $G$, where $G_{i}$ is obtained from $G_{i-1}$
by a single arc adjunction for $i >1$ and $G_{1}$ is isomorphic to an even
subdivision of $O_{1}$, $E_{1}$ or one of the graphs in Figure~\ref{g1}.
Since all possibilities
for $G_{1}$ are either in the list of graphs in Theorem~\ref{main} or
$J'$-compatible
with respect to any assignment $J'$, we may assume that $k >1$. Let 
$H:=G_{k-1}$ and
let $P$ be the unique $\overline{H}$-arc.
Fix a $J$-compatible orientation of $H$ and extend it to an orientation of $G$
arbitrarily. Since $G$ is $J$-incompatible there exist two even 
circuits $A$ and $B$
including $P$ such that $A$ does not have the clockwise parity prescribed by
$J$ but $B$ does. The following lemma shows that the even circuits $A$ and $B$
can be chosen with these properties so that $G[A \cup B]$ is fairly simple.

\medskip

\begin{lem}
\label{AB} The even circuits $A$ and $B$ can be chosen so that
$G[A \cup B]$ is isomorphic to an even subdivision of $O_{1}$,
$O_{2}$, $E_{1}$, $E_{2}$, $E_{3}$, $A_1$, $A_2$, $A_3$, $A_4$ or $A_5$.
(See Figure~\ref{g1}.)
\end{lem}

{\it Proof.} We assume that $A$ and $B$ have been chosen with the
properties above so that $A \cup B$ is minimal.

Let $Q$ be the first $\overline{A} B$-arc we reach when traversing $B$
in a particular direction starting at $P$ and let $R$ be the first
$\overline{A} B$-arc we reach when traversing $B$ in the opposite direction,
again starting at $P$.
If there exists an even circuit in $A \cup Q$ which includes $P \cup Q$ or an
even circuit in $A \cup R$ which includes $P \cup R$ then let $B'$ be this even
circuit. Otherwise there exists an even circuit $B'$ in
$A \cup Q \cup R$ which includes $P \cup Q \cup R$.

First we show that $B'$ has the clockwise parity prescribed by $J$.
Suppose the contrary. The minimality of $A \cup B$ implies
$A \cup B = B' \cup B$ and
therefore $A + B' \subseteq B$. Furthermore $A + B'$ is non-empty and 
the union of circuits.
Therefore $A + B' = B$,  a contradiction to $P \subseteq B$.

If there is a unique $\overline{A} B'$-arc then $G[A \cup B']$ is isomorphic
to an even subdivision of either $O_{1}$ or $E_{1}$. Otherwise
$G[A \cup B']$ is isomorphic to an even subdivision of
$O_{2}$, $E_{2}$, $E_{3}$, $A_1$, $A_2$, $A_3$, $A_4$ or $A_5$.
\qed

\medskip

If $G[A \cup B]$ is an even subdivision of $O_{1}$, $O_{2}$, $E_{1}$,
$E_{2}$ or $E_{3}$ then we have proved Theorem~\ref{main}, for in
these cases $A + B$ is an even circuit with the
clockwise parity prescribed by $J$ since $A+B \subseteq H$. Thus we may
assume that $G[A \cup B]$ is an even subdivision of $A_1$, $A_2$, 
$A_3$, $A_4$ or $A_5$.

\smallskip

\begin{figure}
\begin{center}
\setlength{\unitlength}{1cm}
\hspace{0.5cm}\mbox{\scalebox{0.50}{%
  \includegraphics{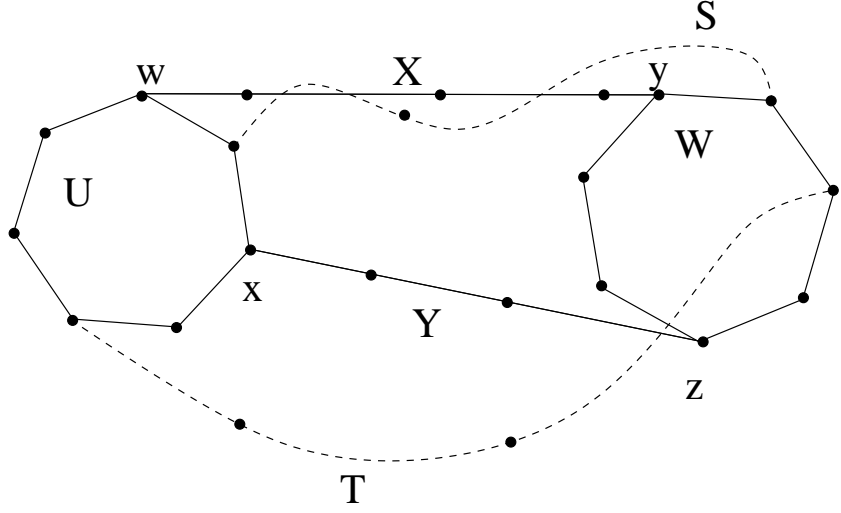}}}%
\end{center}
\caption{$G$ is generated by the even circuits $A$ and $B$ and the two paths
   $S$ and $T$. The paths $S$ and $T$ are dotted because they may 
intersect $X$ and $Y$.}
\label{ABST}
\end{figure}

The symmetric difference $A + B$ is the union of two edge disjoint 
odd circuits $U$ and $W$
in $H$. Since $H$ is even-circuit-connected and therefore 
$2$-connected there exist vertex
disjoint paths $S$ and $T$ in $H$ which join vertices of $U$ to vertices of
$W$ but have no inner vertex in $VU \cup VW$.
Note that $|VU \cap VW| \le 1$ and if $|VU \cap VW|=1$ we may choose
$S$ and $T$ so that
$VS = VU \cap VW$ and $VT \cap VU \cap VW = \emptyset$. Note that
$G[S \cup T \cup U \cup W]$ contains exactly two even circuits $C$
and $D$, and that $C + D = U + W$.
In the following lemma we show that $G$ is spanned by the even circuits
$A$ and $B$ and the paths $S$ and $T$. (See Figure~\ref{ABST}.)

\begin{lem}
\label{key}
$G=G[A \cup B  \cup S \cup T]$.
\end{lem}

{\it Proof.} The set $\{A,B,C,D\}$ of even circuits is
$J$-intractable, for $C$ and $D$ both have the clockwise
parity prescribed by $J$ since $C \cup D \subseteq H$. The assertion
follows by the minimality of $G$. \qed

\medskip

The set $(A \cup B) - (U  \cup W)$ is the union of two vertex disjoint paths
$X$ and $Y$ if $VU \cap VW = \emptyset$. In this case
let $X$ join vertex $w$ in $U$ to vertex $y$ in $W$ and
let $Y$ join vertex $x$ in $U$ to vertex $z$ in $W$. If $VU \cap VW 
\not= \emptyset$ then
let $X$ be the unique path in $(A \cup B) - (U  \cup W)$ and $VY = VU \cap
VW$. Again
we let $X$ join vertex $w$ in $U$ to vertex $y$ in $W$, but we also
write $VU \cap VW = \{x\}$ and $z = x$ in this case.
   (See Figure~\ref{ABST}.)

\begin{figure}
\begin{center}
\setlength{\unitlength}{1cm}
\hspace{0.5cm}\mbox{\scalebox{0.50}{%
  \includegraphics{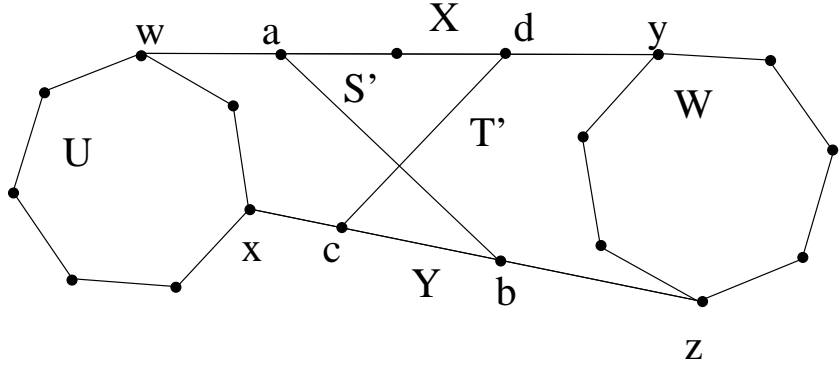}}}%
\end{center}
\caption{Situation in Lemma~\ref{2}.}
\label{ABSITI}
\end{figure}

In our next lemma we show that it is impossible that  the first
$\overline{A \cup B}$-arcs of $S$ and $T$, if we traverse the
paths from their vertex in $U$ to their vertex in $W$, both join
vertices in $VX$ to vertices in $VY$. (See Figure~\ref{ABSITI}.)

\begin{lem}
\label{2} There exist no vertex disjoint $\overline{A \cup
B}$-arcs $S'$ and $T'$ such that $S'$ joins a vertex $a$ in $VX$
to a vertex $b$ in $VY$ and $T'$ joins a vertex $c$ in $VY[x,b] - \{b\}$ to a
vertex $d$ in $VX[a,y]-\{a\}$.
\end{lem}

{\it Proof.}  Clearly the existence of these arcs is impossible if $VU \cap
VW \not= \emptyset$, since $VU \cap VW \not= \emptyset$
implies $VY = \{x\}$.

Suppose therefore that $VU \cap VW = \emptyset$, and that $S'$ and
$T'$ exist. Let $E$ and $F$ be the two even circuits in
   \[ U \cup W \cup X[w,a] \cup X[d,y] \cup Y[x,c] \cup Y[b,z] \cup S' 
\cup  T'. \]
Thus $E + F = U + W = C + D$. If it is not possible to orient $P$ so 
that $E$ and $F$
have the clockwise parity prescribed by $J$ then $\{C,D,E,F\}$ is an
intractable set of circuits. The union of these circuits does not include
$X[a,d] \cup Y[c,b]$, for otherwise $S \cup T$ would include the circuit
$Z = S' \cup T' \cup X[a,d] \cup Y[c,b]$ and this is impossible since 
$S$ and $T$ are
vertex disjoint paths. We now have a contradiction to the minimality of $G$.

Therefore it is
possible to orient $P$ so that $E$ and $F$ have the clockwise parity
prescribed by $J$ and still, by the symmetry of $A$ and $B$, we may 
assume that $A$ does
not have the
prescribed clockwise parity but $B$ does.
Consequently $\{A,B,E,F\}$ is an intractable set of even circuits. By 
the minimality of
$G$ this
implies that $G=G[A \cup B \cup S' \cup T']$.

Suppose that $Z$ is an even circuit. Without loss of
generality we may assume that $A+E=Z$, so that $B+F=Z$. If $Z$ has the
clockwise parity prescribed by $J$ then $\{A,Z,E\}$ is a
$J$-intractable set of even circuits; otherwise $\{B,Z,F\}$ is a
$J$-intractable set of even circuits. This is a contradiction to
the minimality of $G$ for neither $U \cup W \subseteq A \cup E$
nor $U \cup W \subseteq B \cup F$. We conclude that the circuit $Z$ is
odd.

Let $M$ be the unique even circuit that contains $T'$ and is edge
disjoint with $S' \cup W$, let $I$ be the unique even circuit that
contains $T'$ and is edge-disjoint with $S' \cup U$, let $K$ be
the unique even circuit that contains $S'$ and is edge-disjoint
with $T' \cup W$ and, finally, let $L$ be the unique even circuit
that contains $S'$ and is edge-disjoint with $T' \cup U$. Note that
$M + K \ne Z$ since $Z$ is odd. Similarly $I + L \ne Z$. Therefore, since
$M+K+I+L \subseteq U \cup W$  and $M+K+L+I$ is an even cycle, we have
$M+K+L+I = U \cup W$, for if $M + K + L + I = \emptyset$ then
$M + K = I + L = Z$.

Consequently we may assume that $A=M+I$ and $E=M+L$, so that $B=K+L$ and
$F=I+K$. None of $\{A,M,I\}$, $\{B,K,L\}$, $\{E,M,L\}$,
$\{F,I,K\}$ is a $J$-intractable set of even circuits by the
minimality of $G$ for $S' \not\subseteq M \cup I$, $T'
\not\subseteq K \cup L$, $Y[c,b] \not\subseteq M \cup L$ and
$X[a,d] \not\subseteq I \cup K$. Therefore, since $\{A,B,E,F\}$ is
a $J$-intractable set of even circuits but $\{A,M,I\}$ is not,
$\{M,I,B,E,F\}= \{A,B,E,F\} + \{A,M,I\}$ is a $J$-intractable set
of even circuits. It follows that $\{M,I,K,L,E,F\}$ is a $J$-intractable
set of even circuits since $\{B,K,L\}$ is not, and therefore either
$\{E,M,L\}$ or $\{F,I,K\}$
must be a $J$-intractable set of even circuits, a contradiction.
\qed

\medskip

Observe that if $S$, respectively $T$, does not have an
$\overline{A \cup B}$-arc and is therefore equal to either $X$ or
$Y$ then $T$, respectively $S$, has an $\overline{A \cup B}$-arc.
Since $S$ and $T$ are vertex disjoint this arc does not join a
vertex in $VX$ to vertex in $VY$. This fact and Lemma~\ref{2} imply
that, if we traverse
$S$ and $T$ from their vertex of $U$ to their vertex of $W$, then
either the first $\overline{A \cup B}$-arc of $S$ or the first 
$\overline{A \cup B}$-arc of
$T$ does not join a vertex in $VX$ to a vertex in $VY$. Let $R$ denote this arc
for the rest of the paper.
The next lemma shows that one end of $R$ is in $VX \cup VY$.

\begin{lem}
There is no $\overline{A \cup B}$-arc that joins a vertex in $VU - \{w,x\}$
to a vertex in $VW - \{y,z\}$.
\end{lem}

{\it Proof.}
First consider the case that $VU \cap VW = \emptyset$.
Let $Q$ be an $\overline{A \cup B}$-arc that joins a vertex
in $VU - \{w,x\}$ to a vertex in $VW - \{y,z\}$.

Let $E$ and $F$ be the two even circuits in
$U \cup W \cup Q \cup Y$, where we assume without loss of generality that
$P \subseteq X$.  Since $E \cup F \subseteq H$, $E$ and $F$ are of
the prescribed clockwise parity and $\{A , B , E , F\}$ is a
$J$-intractable set of even circuits. Therefore $G=G[A \cup B \cup 
Q]$. Observe that if
$G[A \cup E]=G$ or
$G[B \cup F] = G$ then $G[A \cup F] \not= G$ and $G[B \cup E] \not= G$. By the
symmetry of $E$ and $F$ we therefore assume that $G[A \cup E] \not= G$ and
$G[B \cup F] \not= G$.

The symmetric difference $A + E$ is an even circuit in $U \cup W \cup 
Q \cup X$.
Depending on the clockwise parity of $A+E$ either $\{A,E,A+E\}$ or
$\{A+E,B,F\}$ is a $J$-intractable set of even circuits and therefore
either $G=G[A \cup E]$ or $G=G[B \cup F]$, a contradiction.

Now we consider the case that $VU \cap VW \not= \emptyset$ and, therefore,
$x=z$ and $P \subseteq X$.
Let $E$ and $F$ be the two even circuits in $U \cup W \cup Q$. Since
$E \cup F \subseteq H$,  $E$ and $F$ have the clockwise parity 
prescribed by $J$ and thus $\{A,B,E,F\}$ is a $J$-intractable set of circuits.
There exists at least one even circuit $M$
that includes $Q$ and $X$. Moreover either $A+E=M$ or $A+F=M$. Without loss of
generality let $A+E=M$. Then either $\{A,E,M\}$
or $\{B,F,M\}$ is a $J$-intractable set of circuits, which is a 
contradiction to the
minimality of $G$ for $U \cup W \not\subseteq A \cup E$ and $U \cup W 
\not\subseteq B \cup F$.
\qed

\medskip

In the following lemma we show that both ends of $R$ are either in
$VX$ or in $VY$.

\begin{figure}
\begin{center}
\setlength{\unitlength}{1cm}
\hspace{0.5cm}\mbox{\scalebox{0.50}{%
  \includegraphics{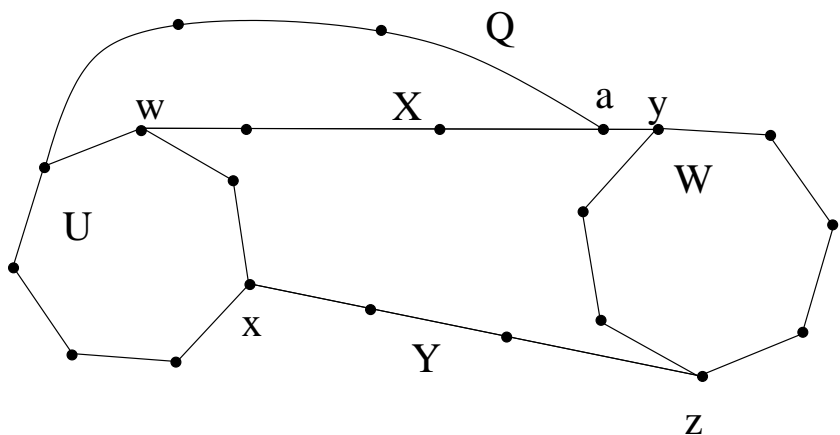}}}%
\end{center}
\caption{Situation in Lemma~\ref{f1}.}
\label{1}
\end{figure}

\begin{lem}
\label{f1}
There is no  $\overline{A \cup B}$-arc that joins a vertex in $(VU
\cup VW) - \{w,x,y,z\}$ to a vertex in $VX \cup VY$.
\end{lem}

{\it Proof.}
Without loss of generality we assume that $Q$ is an  $\overline{(A
\cup B)}$-arc that joins a vertex in $VU - \{w,x\} $ to a vertex $a$ in
$VX - \{w\}$. (See Figure~\ref{1}.)

Let $E$ and $F$ be the two even circuits in $U \cup W \cup Q \cup X[a,y]
\cup Y$. If it is not possible to orient $P$ so that $E$ and $F$ have the
clockwise parity prescribed by $J$ then $\{C,D,E,F\}$
would be a $J$-intractable set of circuits. This conclusion would be 
a contradiction to the
minimality of $G$: $X[w,a]$ is not contained in $C
\cup D \cup E \cup F$ since it cannot be contained in $S \cup T$
because $S$ and $T$ are vertex disjoint.

Therefore we may assume that $E$ and $F$ have the prescribed clockwise parity
and $\{A,B,E,F\}$ is a $J$-intractable set of even circuits. Either 
$A+E$ or $A+F$
is equal to the unique even circuit in $U \cup Q \cup X$. Without loss of
generality we assume that $A+E$ is equal to this even circuit. Then either
$\{A,E,A+E\}$ or $\{A+E, B,F\}$ is a $J$-intractable set of even 
circuits, which is a
contradiction to the minimality of $G$ for neither $W \subseteq A \cup E$ nor
$W \subseteq B \cup F$. \qed

\medskip

\begin{figure}
\begin{center}
\setlength{\unitlength}{1cm}
\hspace{0.5cm}\mbox{\scalebox{0.50}{%
  \includegraphics{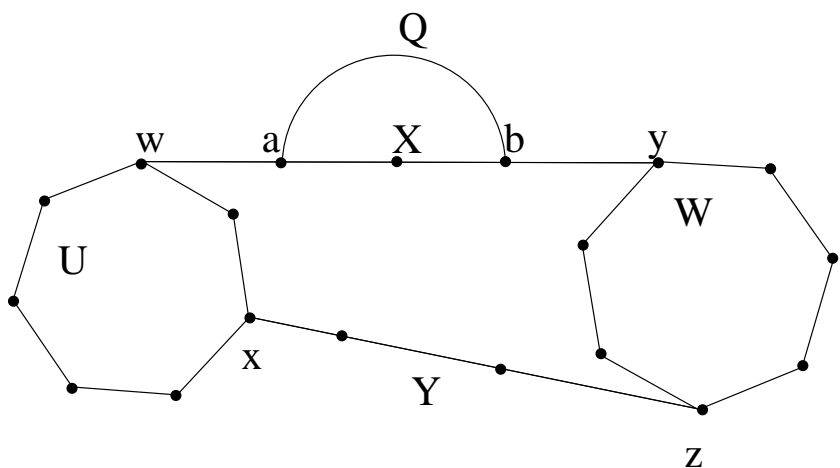}}}%
\end{center}
\caption{Situation in Lemma~\ref{final}}
\label{ABQ}
\end{figure}

Thus $R$ joins either two vertices in $VX$ or two vertices in
$VY$, as in Figure~\ref{ABQ}. In the following lemma we show that $G$ 
is generated by the
even circuits $A$ and $B$ and by the arc $R$ and that $G$ is an even
splitting of $\Delta_1$, $\Delta_2$, $\Delta_3$ or $\Delta_4$.

\begin{lem}
\label{final}
If there is an $\overline{A \cup B}$-arc that joins two vertices in
$VX$ or two vertices in $VY$, then $G$ is an even splitting of
$\Delta_1$, $\Delta_2$, $\Delta_3$ or $\Delta_4$
and $J$ is an assignment which prescribes an odd number of
clockwise even circuits to the even circuits of $G$.
\end{lem}

{\it Proof.} Let $Q$ be an
$\overline{A \cup B}$-arc which joins two
vertices $a$ and $b$ in $VX$. We
assume that $a \in VX[w,b]$. (See Figure~\ref{ABQ}.)

First we show that $G=G[A \cup B \cup Q]$. Let $E$ and $F$ be the
two even circuits in $U \cup W \cup X[w,a] \cup Q \cup X[b,y] \cup
Y $. As in the proofs of the previous lemmas we may orient $P$ such that
$E$ and $F$ have the prescribed clockwise parity. Otherwise $G=G[C 
\cup D \cup E \cup
F]$ and we reach a contradiction: the fact that $S$ and $T$ are
vertex disjoint implies that $Q \cup X[a,b] \not\subseteq S \cup T$, so
that $X[a,b] \not\subseteq C \cup D \cup E \cup F$.
We conclude that $\{A,B,E,F\}$ is a $J$-intractable set of even circuits.

Suppose $Q \cup X[a,b]$ is an even circuit $M$.
Then either $A+E=M$ or $A+F=M$, and without loss of generality we
assume $A+E=M$. Consequently, either $\{A,E,M\}$ or $\{B,F,M\}$ is
a $J$-intractable set of even circuits. We now have a contradiction since
neither $U \cup W \subseteq A \cup E$ nor $U \cup W \subseteq B
\cup F$.

Thus $Q \cup X[a,b]$ is odd and $G$ is an even splitting of
$\Delta_1$, $\Delta_2$, $\Delta_3$ or $\Delta_4$.
\qed

\end{document}